\documentclass[11pt]{article}

\usepackage{amsfonts,amsmath,latexsym,color,epsfig,hyperref,tikz}
\usetikzlibrary {positioning}
\newtheorem{theorem}{Theorem}[section]
\newtheorem{lemma}{Lemma}[section]

\newtheorem{conjecture}{Conjecture}[section]

\newtheorem{definition}{Definition}[section]

\newenvironment{proof}
      {\medskip\noindent{\bf Proof:}\hspace{1mm}}
      {\hfill$\Box$\medskip}
\def\qed{\ifvmode\mbox{ }\else\unskip\fi\hskip 1em plus 10fill$\Box$}

\input{epsf}

\makeatletter
\def\Ddots{\mathinner{\mkern1mu\raise\p@
\vbox{\kern7\p@\hbox{.}}\mkern2mu
\raise4\p@\hbox{.}\mkern2mu\raise7\p@\hbox{.}\mkern1mu}}
\makeatother

\def\mH{\mathcal{H}}

\newcommand*{\eqdef}{\stackrel{\text{\tiny{def}}}{=}}            

\title{\vspace{-0.7cm}Rational exponents in extremal graph theory}
\author{Boris Bukh\thanks{Department of Mathematical Sciences, Carnegie Mellon University, Pittsburgh, PA 15213, USA. E-mail: {\tt bbukh@math.cmu.edu}. Research supported in part by a Sloan Research Fellowship, NSF grant DMS-1301548, and NSF CAREER grant DMS-1555149.}\and David Conlon\thanks{Mathematical Institute, Oxford OX2 6GG, United Kingdom. E-mail: {\tt david.conlon@maths.ox.ac.uk}. Research supported by a Royal Society University Research Fellowship and ERC Starting Grant 676632.}}
\date{}

\begin{document}
\maketitle

\begin{abstract}
Given a family of graphs $\mathcal{H}$, the extremal number $\textrm{ex}(n, \mathcal{H})$ is the largest $m$ for which there exists a graph with $n$ vertices and $m$ edges containing no graph from the family $\mathcal{H}$ as a subgraph. We show that for every rational number $r$ between $1$ and $2$, there is a family of graphs $\mH_r$ such that $\textrm{ex}(n, \mH_r) = \Theta(n^r)$. This solves a longstanding problem in the area of extremal graph theory.
\end{abstract}

\section{Introduction}

Given a family of graphs $\mathcal{H}$, another graph $G$ is said to be $\mathcal{H}$-free if it contains no graph from the family $\mathcal{H}$ as a subgraph. The extremal number $\textrm{ex}(n, \mathcal{H})$ is then defined to be the largest number of edges in an $\mathcal{H}$-free graph on $n$ vertices. If $\mathcal{H}$ consists of a single graph $H$, the classical Erd\H{o}s--Stone--Simonovits theorem~\cite{ES66, ES46} gives a satisfactory first estimate for this function, showing that
\[\textrm{ex}(n, H) = \left(1 - \frac{1}{\chi(H)-1} + o(1)\right) \binom{n}{2},\]
where $\chi(H)$ is the chromatic number of $H$. 

When $H$ is bipartite, the estimate above shows that $\textrm{ex}(n, H) = o(n^2)$. This bound is easily improved to show that for every bipartite graph $H$ there is some positive $\delta$ such that $\textrm{ex}(n, H) = O(n^{2-\delta})$. However, there are very few bipartite graphs for which we have matching upper and lower bounds. 

The most closely studied case is when $H = K_{s,t}$, the complete bipartite graph with parts of order $s$ and $t$. In this case, a famous result of  K\H{o}v\'ari, S\'os and Tur\'an \cite{KST54} shows that $\textrm{ex}(n, K_{s,t}) = O_{s,t}(n^{2 - 1/s})$ whenever $s \leq t$. This bound was shown to be tight for $s = 2$ by Esther Klein~\cite{E38} (see also~\cite{B66, ERS66}) and for $s = 3$ by Brown~\cite{B66}. For higher values of $s$, it is only known that the bound is tight when $t$ is sufficiently large in terms of $s$. This was first shown by Koll\'ar, R\'onyai and Szab\'o \cite{KRS96}, though their construction was improved slightly by Alon, R\'onyai and Szab\'o \cite{ARS99}, who showed that there are graphs with $n$ vertices and $\Omega_s(n^{2 - 1/s})$ edges containing no copy of $K_{s, t}$ with $t = (s-1)! + 1$. 

Alternative proofs showing that $\textrm{ex}(n, K_{s,t}) = \Omega_s(n^{2 - 1/s})$ when $t$ is significantly larger than $s$ were later found by Blagojevi\'c, Bukh and Karasev~\cite{BBK13} and by Bukh~\cite{B14}. In both cases, the basic idea behind the construction is to take a random polynomial $f : \mathbb{F}_q^s \times \mathbb{F}_q^s \rightarrow \mathbb{F}_q$ and then to consider the graph $G$ between two copies of $\mathbb{F}_q^s$ whose edges are all those pairs $(x, y)$ such that $f(x, y) = 0$. A further application of this random algebraic technique was recently given by Conlon~\cite{C14}, who showed that for every natural number $k \geq 2$ there exists a natural number $\ell$ such that, for every $n$, there is a graph on $n$ vertices with $\Omega_k(n^{1 + 1/k})$ edges for which there are at most $\ell$ paths of length $k$ between any two vertices. By a result of Faudree and Simonovits~\cite{FS83}, this is sharp up to the implied constant. We refer the interested reader to~\cite{C14} for further background and details.

In this paper, we give yet another application of the random algebraic method, proving that for every rational number between $1$ and $2$, there is a family of graphs $\mathcal{H}_r$ for which $\textrm{ex}(n, \mathcal{H}_r) = \Theta(n^r)$. This solves a longstanding open problem in extremal graph theory that has been reiterated by a number of authors, including Frankl~\cite{F86} and F\"uredi and Simonovits~\cite{FS13}.

\begin{theorem} \label{main}
For every rational number $r$ between $1$ and $2$, there exists a family of graphs $\mathcal{H}_r$ such that $\textrm{ex}(n, \mathcal{H}_r) = \Theta(n^r)$.
\end{theorem}

Prior to our work, the main result in this direction was due to Frankl~\cite{F86}, who showed that for any rational number $r \geq 1$ there exists a family of $k$-uniform hypergraphs whose extremal function is $\Theta(n^r)$. However, in Frankl's work, the uniformity $k$ depends on the desired exponent $r$, whereas we can always take $k = 2$. 



In order to define the relevant families $\mathcal{H}_r$, we need some preliminary definitions.

\begin{definition}
A rooted tree $(T, R)$ consists of a tree $T$ together with an independent set $R \subset V(T)$, which we refer to as the roots. When the set of roots is understood, we will simply write $T$.
\end{definition}

Each of our families $\mathcal{H}_r$ will be of the following form.

\begin{definition}
Given a rooted tree $(T, R)$, we define the $p$th power $\mathcal{T}_R^p$ of $(T, R)$ to be the family of graphs consisting of all possible unions of $p$ distinct labelled copies of $T$, each of which agree on the set of roots $R$. Again, we will usually omit $R$, denoting the family by $\mathcal{T}^p$ and referring to it as the $p$th power of $T$.
\end{definition}

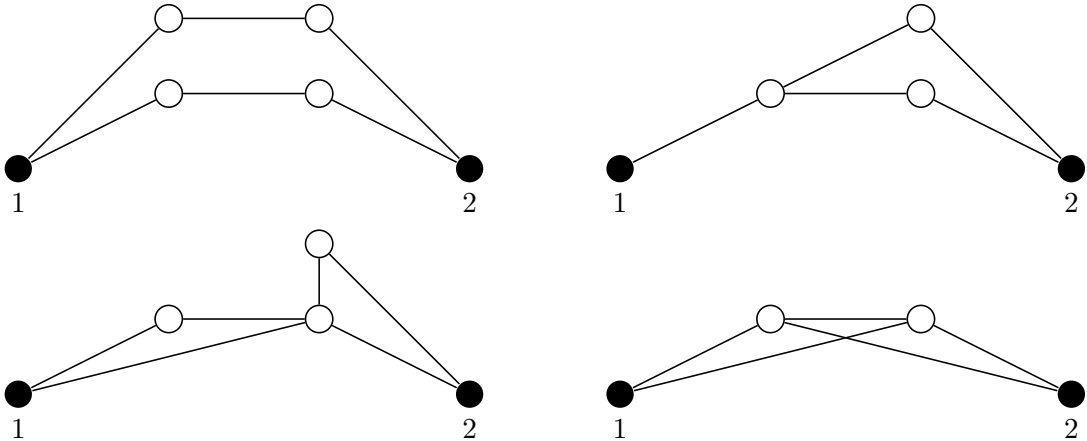
\begin{figure}[t!]
\definecolor {processblue}{cmyk}{0.96,0,0,0}
\begin {center}
\begin {tikzpicture}[-latex, auto, node distance = 1 cm and 2cm, on grid,
semithick, 
state/.style ={ circle, top color =white, 
draw, black, text=blue, minimum width = 1 mm}]

\begin{scope}[shift={(11,0)}]
\node[state] (B) {};
\node[circle, fill = black] (A) [below left=of B,label=below:$1$] {};
\node[state] (C) [right =of B] {};
\node[circle, fill = black] (D) [below right = of C,label=below:$2$] {};
\node[state] (E) [above = of B] {};
\node[state] (F) [above = of C] {};

\path[-] (A) edge (B);
\path[-] (A) edge (E);
\path[-] (B) edge (C);
\path[-] (C) edge (D);
\path[-] (E) edge (F);
\path[-] (D) edge (F);
\end{scope}

\begin{scope}[shift={(19,0)}]
\node[state] (B) {};
\node[circle, fill = black] (A) [below left=of B,label=below:$1$] {};
\node[state] (C) [right =of B] {};
\node[circle, fill = black] (D) [below right = of C,label=below:$2$] {};
\node[state] (F) [above = of C] {};

\path[-] (A) edge (B);
\path[-] (B) edge (C);
\path[-] (C) edge (D);
\path[-] (B) edge (F);
\path[-] (D) edge (F);
\end{scope}

\begin{scope}[shift={(11,-3)}]
\node[state] (B) {};
\node[circle, fill = black] (A) [below left=of B,label=below:$1$] {};
\node[state] (C) [right =of B] {};
\node[circle, fill = black] (D) [below right = of C,label=below:$2$] {};
\node[state] (F) [above = of C] {};

\path[-] (A) edge (B);
\path[-] (A) edge (C);
\path[-] (B) edge (C);
\path[-] (C) edge (D);
\path[-] (C) edge (F);
\path[-] (D) edge (F);
\end{scope}

\begin{scope}[shift={(19,-3)}]
\node[state] (B) {};
\node[circle, fill = black] (A) [below left=of B,label=below:$1$] {};
\node[state] (C) [right =of B] {};
\node[circle, fill = black] (D) [below right = of C,label=below:$2$] {};

\path[-] (A) edge (B);
\path[-] (A) edge (C);
\path[-] (B) edge (C);
\path[-] (C) edge (D);
\path[-] (B) edge (D);
\end{scope}

\end{tikzpicture}
\end{center}
\caption{Some of the graphs in $\mathcal{T}^2$ when $(T, R)$ is a path of length $3$ with rooted endpoints. The remaining graphs in $\mathcal{T}^2$ are obtained
by swapping the two roots, which are labelled $1$ and $2$.} \label{fig:T2}
\end{figure}

We note that $\mathcal{T}^p$ consists of more than one graph because we allow the unrooted vertices $V(T)\setminus R$ to meet in every possible way. For example, if $T$ is a path of length $3$ whose endpoints are rooted, the family $\mathcal{T}^2$ contains a cycle of length $6$ and the various degenerate configurations shown in Figure~\ref{fig:T2}.

The following parameter will be critical in studying the extremal number of the family $\mathcal{T}^p$.

\begin{definition}
Given a rooted tree $(T, R)$, we define the density $\rho_T$ of $(T, R)$ to be $\frac{e(T)}{v(T) - |R|}$.
\end{definition}

The upper bound in Theorem~\ref{main} will follow from an application of the next lemma.

\begin{lemma} \label{upper}
For any rooted tree $(T, R)$ with at least one root, the family $\mathcal{T}^p$ satisfies
\[\textrm{ex}(n, \mathcal{T}^p) = O_p(n^{2 - 1/\rho_T}).\]
\end{lemma}

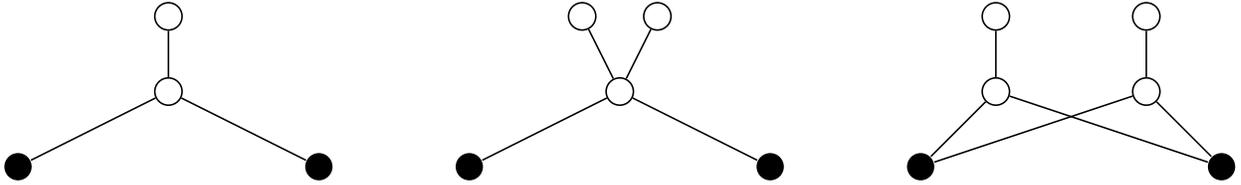
\begin{figure}
\definecolor {processblue}{cmyk}{0.96,0,0,0}
\begin {center}
\begin {tikzpicture}[-latex, auto, node distance = 1 cm and 2cm, 
on grid, 
semithick, 
state/.style ={ circle, top color =white, draw, black, text=blue, minimum width = 1 mm}]

\begin{scope}[shift={(9,0)}]
\node[state] (B) {};
\node[circle, fill = black] (A) [below left=of B] {};
\node[circle, fill = black] (C) [below right =of B] {};
\node[state] (D) [above = of B] {};

\path[-] (A) edge (B);
\path[-] (B) edge (C);
\path[-] (B) edge (D);
\end{scope}

\begin{scope}[shift={(15,0)}]
\node[state] (C) {};
\node[circle, fill = black] (A) [below left=of C] {};
\node[circle, fill = black] (B) [below right=of C] {};
\node[draw, white] (X) [above = of C] {};
\node[state] (D) [left = of X, xshift = 1.5cm] {};
\node[state] (E) [right = of X, xshift = -1.5cm] {};

\path[-] (A) edge (C);
\path[-] (B) edge (C);
\path[-] (D) edge (C);
\path[-] (E) edge (C);
\end{scope}

\begin{scope}[shift={(21,0)}]
\node[draw, white] (X) {};
\node[circle, fill = black] (A) [below left=of X] {};
\node[circle, fill = black] (B) [below right=of X] {};
\node[state] (C) [left = of X, xshift = 1cm] {};
\node[state] (D) [right = of X, xshift = -1cm] {};
\node[state] (E) [above = of C] {};
\node[state] (F) [above = of D] {};

\path[-] (A) edge (C);
\path[-] (B) edge (D);
\path[-] (B) edge (C);
\path[-] (A) edge (D);
\path[-] (C) edge (E);
\path[-] (D) edge (F);
\end{scope}

\end{tikzpicture}
\end{center}
\caption{An unbalanced rooted tree $T$ and two elements of $\mathcal{T}^2$.} \label{fig:unbalanced}
\end{figure}

It would be wonderful if there were also a matching lower bound for $\textrm{ex}(n, \mathcal{T}^p)$. However, this is in general too much to expect. If, for example, $(T, R)$ is the star $K_{1, 3}$ with two rooted leaves, $\mathcal{T}^2$ will contain the graph shown in Figure~\ref{fig:unbalanced} where the two central vertices agree. However, this graph is a tree, so it is easy to show that $\textrm{ex}(n, \mathcal{T}^2) = O(n)$, whereas, since $\rho_T = 3/2$, Lemma~\ref{upper} only gives $\textrm{ex}(n, \mathcal{T}^2) = O(n^{4/3})$. Luckily, we may avoid these difficulties by restricting attention to so-called balanced trees. 

\begin{definition}
Given a subset $S$ of the unrooted vertices $V(T)\setminus R$ in a rooted tree $(T, R)$, we define the density $\rho_S$ of $S$ to be $e(S)/|S|$, where $e(S)$ is the number of edges in $T$ with at least one endpoint in $S$. Note that when $S = V(T)\setminus R$, this agrees with the definition above. We say that the rooted tree $(T, R)$ is balanced if, for every subset $S$ of the unrooted vertices $V(T)\setminus R$, the density of $S$ is at least the density of $T$, that is, $\rho_S \geq \rho_T$. In particular, if $|R|\geq 2$, then this condition guarantees that every leaf in the tree is a root.
\end{definition}

With the caveat that our rooted trees must be balanced, we may now prove a lower bound matching Lemma~\ref{upper} by using the random algebraic method.

\begin{lemma} \label{lower}
For any balanced rooted tree $(T, R)$, there exists a positive integer $p$ such that the family $\mathcal{T}^p$ satisfies
\[\textrm{ex}(n, \mathcal{T}^p) = \Omega(n^{2 - 1/\rho_T}).\]
\end{lemma}

Therefore, given a rational number $r$ between $1$ and $2$, it only remains to identify a balanced rooted tree $(T, R)$ for which $2 - 1/\rho_T$ is equal to $r$.

\begin{definition}
Suppose that $a$ and $b$ are natural numbers satisfying $a-1\leq b<2a-1$ and put $i=b-a$.
We define a rooted tree $T_{a,b}$ by taking a path with $a$ vertices, which are labelled in order as $1, 2, \dots, a$, and then adding an additional rooted leaf to each of the $i+1$ vertices 
\[
  1, \left\lfloor 1+\frac{a}{i}\right\rfloor,\left\lfloor 1+2 \cdot \frac{a}{i}\right\rfloor,\dotsc, \left\lfloor 1+(i-1) \cdot \frac{a}{i}\right\rfloor, a.
\]

For $b\geq 2a-1$, we define $T_{a,b}$ recursively to be the tree obtained by attaching a rooted leaf to each unrooted vertex of~$T_{a,b-a}$.
\end{definition}


Note that the tree $T_{a,b}$ has $a$ unrooted vertices and $b$ edges, so that $\rho_T=b/a$.
Now, given a rational number $r$ with $1 < r <2$, let $a/b = 2 - r$ and let $\mathcal{T}_{a,b}^p$ be the $p$th power of $T_{a,b}$. To prove Theorem~\ref{main}, it will suffice to prove that $T_{a, b}$ is balanced, since we may then apply Lemmas~\ref{upper} and \ref{lower} to $\mathcal{T}_{a,b}^p$, for $p$ sufficiently large, to conclude that 
\[\textrm{ex}(n, \mathcal{T}_{a,b}^p) = \Theta(n^{2 - a/b}) = \Theta(n^r).\]
Therefore, the following lemma completes the proof of Theorem~\ref{main}.

\begin{lemma} \label{balanced}
The tree $T_{a,b}$ is balanced.
\end{lemma}

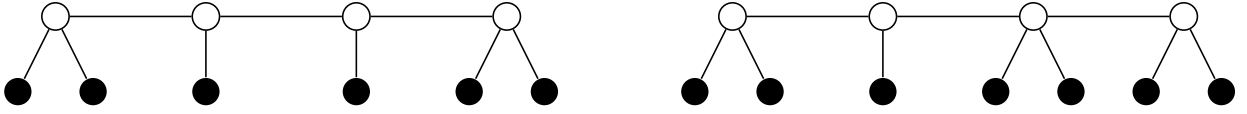
\begin{figure}
\definecolor {processblue}{cmyk}{0.96,0,0,0}
\begin {center}
\begin {tikzpicture}[-latex, auto, node distance = 1 cm and 2cm, 
on grid, 
semithick, 
state/.style ={ circle, top color =white, draw, black, text=blue, minimum width = 1 mm}]

\begin{scope}[shift={(10,0)}]
\node[state] (A) {};
\node[state] (B) [right = of A] {};
\node[state] (C) [right = of B] {};
\node[state] (D) [right = of C] {};
\node[draw = white] (E) [below = of A] {};
\node[circle, fill = black] (F) [below = of B] {};
\node[circle, fill = black] (G) [below = of C] {};
\node[draw = white] (H) [below = of D] {};
\node[circle, fill = black] (I) [left = of E, xshift = 1.5cm] {};
\node[circle, fill = black] (J) [right = of E, xshift = -1.5cm] {};
\node[circle, fill = black] (K) [left = of H, xshift = 1.5cm] {};
\node[circle, fill = black] (L) [right = of H, xshift = -1.5cm] {};

\path[-] (A) edge (B);
\path[-] (B) edge (C);
\path[-] (C) edge (D);
\path[-] (A) edge (I);
\path[-] (A) edge (J);
\path[-] (B) edge (F);
\path[-] (C) edge (G);
\path[-] (D) edge (K);
\path[-] (D) edge (L);
\end{scope}

\begin{scope}[shift={(19,0)}]
\node[state] (A) {};
\node[state] (B) [right = of A] {};
\node[state] (C) [right = of B] {};
\node[state] (D) [right = of C] {};
\node[draw = white] (E) [below = of A] {};
\node[circle, fill = black] (F) [below = of B] {};
\node[draw = white] (G) [below = of C] {};
\node[draw = white] (H) [below = of D] {};
\node[circle, fill = black] (I) [left = of E, xshift = 1.5cm] {};
\node[circle, fill = black] (J) [right = of E, xshift = -1.5cm] {};
\node[circle, fill = black] (K) [left = of H, xshift = 1.5cm] {};
\node[circle, fill = black] (L) [right = of H, xshift = -1.5cm] {};
\node[circle, fill = black] (M) [left = of G, xshift = 1.5cm] {};
\node[circle, fill = black] (N) [right = of G, xshift = -1.5cm] {};

\path[-] (A) edge (B);
\path[-] (B) edge (C);
\path[-] (C) edge (D);
\path[-] (A) edge (I);
\path[-] (A) edge (J);
\path[-] (B) edge (F);
\path[-] (C) edge (M);
\path[-] (C) edge (N);
\path[-] (D) edge (K);
\path[-] (D) edge (L);

\end{scope}

\end{tikzpicture}
\end{center}
\caption{The rooted trees $T_{4,9}$ and $T_{4, 10}$.} \label{fig:balanced}
\end{figure}

All of the proofs will be given in the next section: we will prove the easy Lemma~\ref{upper} in Section~\ref{sec:upper}; Lemma~\ref{balanced} and another useful fact about balanced trees will be proved in Section~\ref{sec:balanced}; and Lemma~\ref{lower} will be proved in Section~\ref{sec:lower}. We conclude, in Section~\ref{sec:conclusion}, with some brief remarks.

\section{Proofs}

\subsection{The upper bound} \label{sec:upper}

We will use the following folklore lemma.

\begin{lemma} \label{mindeg}
A graph $G$ with average degree $d$ has a subgraph $G'$ of minimum degree at least $d/2$.
\end{lemma}

With this mild preliminary, we are ready to prove Lemma~\ref{upper}, that $\textrm{ex}(n, \mathcal{T}^p) = O_p(n^{2 - 1/\rho_T})$ for any rooted tree $(T, R)$.

\vspace{3mm}
{\bf Proof of Lemma~\ref{upper}:} 
Suppose that $G$ is a graph on $n$ vertices with $c n^{2 - \alpha}$ edges, where $\alpha= 1/\rho_T$ and $c \geq 2 \max(|T|, p)$. We wish to show that $G$ contains an element of $\mathcal{T}^p$. Since the average degree of $G$ is $2c n^{1 - \alpha}$, Lemma~\ref{mindeg} implies that $G$ has a subgraph $G'$ with minimum degree at least $c n^{1 - \alpha}$. Suppose that this subgraph has $s \leq n$ vertices. By embedding greedily one vertex at a time, the minimum degree condition allows us to conclude that $G'$ contains at least 
\[s \cdot c n^{1 - \alpha} \cdot (cn^{1- \alpha} - 1) \cdots (cn^{1 - \alpha} - |T| + 2) \geq (c/2)^{|T| - 1} s n^{(|T| - 1)(1 - \alpha)}\]
labelled copies of the (unrooted) tree $T$. Since there are at most $s^{|R|}$ possible choices for the root vertices $R$, there must be some choice $R_0$ for these vertices in at least
\[\frac{(c/2)^{|T| - 1} s n^{(|T| - 1)(1 - \alpha)}}{s^{|R|}} \geq \frac{(c/2)^{|T| - 1} n^{(|T| - 1)(1 - \alpha)}}{n^{|R| - 1}} = (c/2)^{|T| - 1}\]
distinct labelled copies of $T$, where we used that $s \leq n$ and $\alpha = 1/\rho_T = (|T| - |R|)/(|T| - 1)$. Since $(c/2)^{|T| - 1} \geq p$, this gives the required element of $\mathcal{T}^p$.
\qed

\subsection{Balanced trees} \label{sec:balanced}

We will begin by proving Lemma~\ref{balanced}, that $T_{a,b}$ is balanced. 

\vspace{3mm}
{\bf Proof of Lemma~\ref{balanced}:}
Suppose that $S$ is a proper subset of the unrooted vertices of $T_{a,b}$. We wish to show that $e(S)$, the number of edges in $T$ with at least one endpoint in $S$, is at least $\rho_T |S|$,  where $\rho_T = b/a$. We may make two simplifying assumptions. First, we may assume that $a-1\leq b<2a-1$. Indeed, if $b\geq 2a-1$, then the bound for $T_{a,b}$ follows from
the bound for $T_{a,b-a}$, which we may assume by induction. Second, we may assume that the vertices in $S$ form a subpath of the base path of length $a$. Indeed, given the result in this case, we may write any $S$ as the disjoint union of subpaths $S_1, S_2, \dots, S_p$ with no edges between them, so that  
\[e(S) = e(S_1 \cup S_2 \cup \dots \cup S_p) = e(S_1) + e(S_2) + \dots + e(S_p) \geq \rho_T (|S_1| + |S_2| + \dots + |S_p|) = \rho_T |S|.\]
Suppose, therefore, that $S=\{l,l+1,\dotsc,r\}$ is a proper subpath of the base path $\{1, 2, \dots, a\}$ and $b-a = i$. 

As the desired claim is trivially true if $i=-1$, we will assume that $i\geq 0$. In particular, it follows from this assumption that vertex $1$ of the base path is adjacent to a rooted
vertex.

Let $R$ be the number of rooted vertices adjacent to $S$. 
For $0\leq j\leq i - 1$, the $j$th rooted vertex is adjacent to $S$ precisely when $l\leq 1 + j \left( \frac{a}{i} \right)  <r+1$, which is equivalent to
\[
  (l-1)\frac{i}{a}\leq j  <r\frac{i}{a}.
\]
Therefore, if $a$ is not contained in $S$, it follows that $R\geq \lfloor|S|\frac{i}{a}\rfloor=\lfloor|S|\frac{b-a}{a}\rfloor$. Furthermore, if $l = 1$, then $R = \lceil |S|\frac{b-a}{a}\rceil$. Finally, if $r = a$ and $i>0$, then, using 
\[a - \left\lfloor 1 + j \cdot \frac{a}{i} \right\rfloor \leq (i - j)\frac{a}{i},\] 
it follows that $S$ is adjacent to the $j$th root whenever $i|S|/a>i-j$, and so $R \geq \lceil |S|\frac{b-a}{a}\rceil$.

\textit{Case 1:} $i=0$. Since $S$ is a proper subpath, it is adjacent to at least $|S|=(b/a)|S|$ edges.

\textit{Case 2:} $R \geq \lceil|S|\frac{b-a}{a}\rceil$. Then the total number of edges adjacent to $S$ is at least $R+|S|\geq (b/a)|S|$.

\textit{Case 3:} $i>0$ and $R < \lceil|S|\frac{b-a}{a}\rceil$. Then $S$ is adjacent to $|S|+1$ edges in the base path, for a total
of $\lfloor|S|\frac{b-a}{a}\rfloor+|S|+1\geq (b/a)|S|$ adjacent edges.
\qed

\vspace{3mm}
Before moving on to the proof of Lemma~\ref{lower}, it will be useful to note that if $T$ is balanced then every graph in $\mathcal{T}^p$ is at least as dense as $T$.

\begin{lemma} \label{Hbalanced}
If $(T, R)$ is a balanced rooted tree, then every graph $H$ in $\mathcal{T}^s$ satisfies 
\[e(H) \geq \rho_T(|H| - |R|).\]
\end{lemma}

\begin{proof}
We will prove the result by induction on $s$. It is clearly true when $s = 1$, so we will assume that it holds for any $H \in \mathcal{T}^s$ and prove it when $H \in \mathcal{T}^{s+1}$.

Suppose, therefore, that $H$ is the union of $s+1$ labelled copies of $T$, say $T_1, T_2, \dots, T_{s+1}$, each of which agree on the set of roots $R$. If we let $H'$ be the union of the first $s$ copies of $T$, the induction hypothesis tells us that $e(H') \geq \rho_T(|H'| - |R|)$. Let $S$ be the set of vertices in $T_{s+1}$ which are not contained in $H'$. Then, since $T$ is balanced, we know that $e(S)$, the number of edges in $T_{s+1}$ (and, therefore, in $H$) with at least one endpoint in $S$, is at least $\rho_T |S|$. It follows that
\[e(H) \geq e(H') + e(S) \geq \rho_T (|H'| - |R|) + \rho_T |S| = \rho_T (|H| - |R|),\]
as required.
\end{proof}

\subsection{The lower bound} \label{sec:lower}

The proof of the lower bound will follow~\cite{B14} and \cite{C14} quite closely. We begin by describing the basic setup and stating a number of lemmas which we will require in the proof. We will omit the proofs of these lemmas, referring the reader instead to~\cite{B14} and~\cite{C14}.

Let $q$ be a prime power and let $\mathbb{F}_q$ be the finite field of order $q$. We will consider polynomials in $t$ variables over $\mathbb{F}_q$, writing any such polynomial as $f(X)$, where $X = (X_1, \dots, X_t)$. We let $\mathcal{P}_d$ be the set of polynomials in $X$ of degree at most $d$, that is, the set of linear combinations over $\mathbb{F}_q$ of monomials of the form $X_1^{a_1} \cdots X_t^{a_t}$ with $\sum_{i=1}^t a_i \leq d$. By a random polynomial, we just mean a polynomial chosen uniformly from the set $\mathcal{P}_d$. One may produce such a random polynomial by choosing the coefficients of the monomials above to be random elements of $\mathbb{F}_q$.

The first result we will need says that once $q$ and $d$ are sufficiently large, the probability that a randomly chosen polynomial from $\mathcal{P}_d$ contains each of $m$ distinct points is exactly $1/q^{m}$.

\begin{lemma} \label{graphexp}
Suppose that $q > \binom{m}{2}$ and $d \geq m - 1$. Then, if $f$ is a random polynomial from $\mathcal{P}_d$ and $x_1, \dots, x_m$ are $m$ distinct points in $\mathbb{F}_q^t$,
$$\mathbb{P}[f(x_i) = 0 \mbox{ for all } i = 1, \dots, m] = 1/q^{m}.$$ 
\end{lemma}

We also need to note some basic facts about affine varieties over finite fields. If we write $\overline{\mathbb{F}}_q$ for the algebraic closure of $\mathbb{F}_q$, a variety over $\overline{\mathbb{F}}_q$ is a set of the form
\[W = \{x \in \overline{\mathbb{F}}_q^t : f_1(x) = \dots = f_s(x) = 0\}\]
for some collection of polynomials $f_1, \dots, f_s \colon \overline{\mathbb{F}}_q^t \rightarrow \overline{\mathbb{F}}_q$. We say that $W$ is defined over $\mathbb{F}_q$ if the coefficients of these polynomials are in $\mathbb{F}_q$ and write $W(\mathbb{F}_q) = W \cap \mathbb{F}_q^t$. We say that $W$ has complexity at most $M$ if $s$, $t$ and the degrees of the $f_i$ are all bounded by $M$. Finally, we say that a variety is absolutely irreducible if it is irreducible over $\overline{\mathbb{F}}_q$, reserving the term irreducible for irreducibility over $\mathbb{F}_q$ of varieties defined over $\mathbb{F}_q$.

The next result we will need is the Lang--Weil bound \cite{LW54} relating the dimension of a variety $W$ to the number of points in $W(\mathbb{F}_q)$. It will not be necessary to give a formal definition for the dimension of a variety, though some intuition may be gained by noting that if $f_1, \dots, f_s\colon \overline{\mathbb{F}}_q^t \rightarrow \overline{\mathbb{F}}_q$ are generic polynomials then the dimension of the variety they define is $t - s$. 

\begin{lemma} \label{LW}
Suppose that $W$ is a variety over $\overline{\mathbb{F}}_q$ of complexity at most $M$. Then 
$$|W(\mathbb{F}_q)| = O_M(q^{\dim W}).$$
Moreover, if $W$ is defined over $\mathbb{F}_q$ and absolutely irreducible, then 
$$|W(\mathbb{F}_q)| = q^{\dim W}(1 + O_M(q^{-1/2})).$$
\end{lemma}

We will also need the following standard result from algebraic geometry, which says that if $W$ is an absolutely irreducible variety and $D$ is a variety intersecting $W$, then either $W$ is contained in $D$ or its intersection with $D$ has smaller dimension.

\begin{lemma} \label{dimdrop}
Suppose that $W$ is an absolutely irreducible variety over $\overline{\mathbb{F}}_q$ and $\dim W \geq 1$. Then, for any variety $D$, either $W\subseteq D$ or $W\cap D$ is a variety
of dimension less than $\dim W$. 
\end{lemma}

The final ingredient we require says that if $W$ is a variety which is defined over $\mathbb{F}_q$, then there is a bounded collection of absolutely irreducible varieties $Y_1, \dots, Y_t$, each of which is defined over $\mathbb{F}_q$, such that $\cup_{i=1}^t Y_i(\mathbb{F}_q) = W(\mathbb{F}_q)$. 

\begin{lemma} \label{frob}
Suppose that $W$ is a variety over $\overline{\mathbb{F}}_q$ of complexity at most $M$ which is defined over $\mathbb{F}_q$. Then there are $O_M(1)$ absolutely irreducible varieties $Y_1, \dots, Y_t$, each of which is defined over $\mathbb{F}_q$ and has complexity $O_M(1)$, such that $\cup_{i=1}^t Y_i(\mathbb{F}_q) = W(\mathbb{F}_q)$.
\end{lemma}

We can combine the preceding three lemmas into a single result as follows:
\begin{lemma}\label{dichotomy}
Suppose $W$ and $D$ are varieties over $\overline{\mathbb{F}}_q$ of complexity at most $M$ which are defined over $\mathbb{F}_q$.
Then one of the following holds for all $q$ sufficiently large in terms of $M$:
\begin{itemize}
\item $|W(\mathbb{F}_q)\setminus D(\mathbb{F}_q)|\geq q/2$, or
\item $|W(\mathbb{F}_q)\setminus D(\mathbb{F}_q)|\leq c$, where $c=c_M$ depends only on $M$.
\end{itemize}
\end{lemma}

\begin{proof}
By Lemma~\ref{frob}, there is a decomposition $W(\mathbb{F}_q)=\bigcup_{i=1}^t Y_i(\mathbb{F}_q)$ for some
bounded-complexity absolutely irreducible varieties $Y_i$ defined over $\mathbb{F}_q$. If $\dim Y_i\geq 1$,
Lemma~\ref{dimdrop} tells us that either $Y_i\subset D$ or the dimension of $Y_i\cap D$ is smaller than the dimension of $Y_i$.
If $Y_i\subset D$, then the component does not contribute any point to $W(\mathbb{F}_q)\setminus D(\mathbb{F}_q)$
and may be discarded. If instead the dimension of $Y_i\cap D$ is smaller than the dimension of $Y_i$, the Lang--Weil bound, Lemma~\ref{LW}, tells us that for $q$ sufficiently
large
\[
  |W(\mathbb{F}_q)\setminus D(\mathbb{F}_q)|\geq |Y_i(\mathbb{F}_q)|-|Y_i(\mathbb{F}_q)\cap D|\geq q^{\dim Y_i}-O(q^{\dim Y_i-\tfrac{1}{2}})-O(q^{\dim Y_i-1})\geq q/2.
\]

On the other hand, if $\dim Y_i=0$ for every $Y_i$ which is not contained in $D$, Lemma~\ref{LW} tells us that $|W(\mathbb{F}_q)\setminus D(\mathbb{F}_q)|\leq \sum |Y_i(\mathbb{F}_q)| =O(1)$, where the sum
is taken over all $i$ for which $\dim Y_i=0$.
\end{proof}

We are now ready to prove Lemma~\ref{lower}, that for any balanced rooted tree $(T, R)$ there exists a positive integer $p$ such that $\textrm{ex}(n, \mathcal{T}^p) = \Omega(n^{2 - 1/\rho_T})$.

\vspace{3mm}
{\bf Proof of Lemma~\ref{lower}:}
Let $(T, R)$ be a balanced rooted tree with $a$ unrooted vertices and $b$ edges, where $R = \{u_1, \dots, u_r\}$ and $V(T) \setminus R = \{v_1, \dots, v_a\}$. Let $s = 2br$, $d = s b$, $N = q^b$ and suppose that $q$ is sufficiently large. Let $f_1, \dots, f_{a} \colon \mathbb{F}_q^b \times \mathbb{F}_q^b \rightarrow \mathbb{F}_q$ be independent random polynomials in $\mathcal{P}_d$. We will consider the bipartite graph $G$ between two copies $U$ and $V$ of $\mathbb{F}_q^b$, each of order $N = q^b$, where $(u, v)$ is an edge of $G$ if and only if 
$$f_1(u, v) = \dots = f_{a}(u,v) = 0.$$ 
Since $f_1, \dots, f_{a}$ were chosen independently, Lemma~\ref{graphexp} with $m = 1$ tells us that the probability a given edge $(u,v)$ is in $G$ is $q^{-a}$. Therefore, the expected number of edges in $G$ is $q^{-a} N^2 = N^{2 - a/b}$.

Suppose now that $w_1, w_2, \dots, w_r$ are fixed vertices in $G$ and let $C$ be the collection of copies of $T$ in $G$ such that $w_i$ corresponds to $u_i$ for all $1 \leq i \leq r$. We will be interested in estimating the $s$-th moment of $|C|$. To begin, we note that $|C|^s$ counts the number of ordered collections of $s$ (possibly overlapping or identical) copies of $T$ in $G$ such that $w_i$ corresponds to $u_i$ for all $1 \leq i \leq r$. Since the total number of edges $m$ in a given collection of $s$ rooted copies of $T$ is at most $s b$ and $q$ is sufficiently large, Lemma~\ref{graphexp} tells us that the probability this particular collection of copies of $T$ is in $G$ is $q^{-a m}$, where we again use the fact that $f_1, \dots, f_a$ are chosen independently. 

Suppose that $H$ is an element of $\mathcal{T}^{s}_{\leq} \eqdef \mathcal{T}^1 \cup \mathcal{T}^2 \cup \dots \cup \mathcal{T}^s$. Within the complete bipartite graph from $U$ to $V$, let $N_s(H)$ be the number of ordered collections of $s$ copies of $T$, each rooted at $w_1, \dots, w_r$ in the same way, whose union is a copy of $H$. Then
\[\mathbb{E}[|C|^s] = \sum_{H \in \mathcal{T}^{s}_{\leq}} N_s(H) q^{-a e(H)},\]
while $N_s(H) = O_s(N^{|H| - |R|})$. Since $(T, R)$ is balanced, Lemma~\ref{Hbalanced} shows that $\frac{e(H)}{|H| - |R|} \geq \rho_T = \frac{b}{a}$ for every $H \in \mathcal{T}^s_{\leq}$. It follows that
\begin{eqnarray*}
\mathbb{E}[|C|^s] & = & \sum_{H \in \mathcal{T}^s_{\leq}} N_s(H) q^{-a e(H)}  =  \sum_{H \in \mathcal{T}^s_{\leq}} O_s\left(N^{|H| - |R|}\right) q^{-a e(H)}\\
& = & O_s \left(\sum_{H \in \mathcal{T}^s_{\leq}} q^{b(|H|-|R|)} q^{-a e(H)} \right) = O_s(1).
\end{eqnarray*}
By Markov's inequality, we may conclude that
\[\mathbb{P}[|C| \geq c] = \mathbb{P}[|C|^s \geq c^s] \leq \frac{\mathbb{E}[|C|^s]}{c^s} = \frac{O_s(1)}{c^s}.\]
Our aim now is to show that $|C|$ is either quite small or very large. To begin, note that the set $C$ is a subset of $X(\mathbb{F}_q)$, where $X$ is the algebraic variety defined as the set of $(x_1, \dots, x_a) \in \overline{\mathbb{F}}_q^{ba}$ satisfying the equations
\begin{itemize}
\item
$f_i(w_k, x_\ell) = 0$ for all $k$ and $\ell$ such that $(u_k, v_{\ell}) \in T$ and
\item
$f_i(x_k, x_\ell) = 0$ for all $k$ and $\ell$ such that $(v_k, v_{\ell}) \in T$
\end{itemize}
for all $i = 1, 2, \dots, a$. For each $i \neq j$ such that $v_i$ and $v_j$ are on the same side of the natural bipartition of $T$, we let
\[D_{ij} = X \cap \{(x_1, \dots, x_a) : x_i = x_j\}\]
and, for each $k, \ell$ such that $v_k$ and $u_\ell$ are on the same side of the bipartition, we let
\[D'_{k\ell} = X \cap \{(x_1, \dots, x_a) : x_k = w_\ell\}.\]
We put
\[
  D\eqdef\bigcup_{i, j} D_{ij} \cup \bigcup_{k,\ell} D'_{k\ell}.
\]
The sets $D_{ij}$ and $D'_{k \ell}$ capture those elements of $X$ which are degenerate and so not elements of $C$. As a union of varieties
is a variety, the set $D$ is a variety that captures all degenerate elements of~$X$. Furthermore, the complexity
of $D$ is bounded since the number and complexity of the $D_{ij}$ and $D'_{k\ell}$ is bounded.

By Lemma~\ref{dichotomy}, we see that that there exists a constant $c_T$, depending only on $T$, such that either $|C| \leq c_T$ or $|C| \geq q/2$. Therefore, by the consequence of Markov's inequality noted earlier,
\[\mathbb{P}[|C| > c_T] = \mathbb{P}[|C| \geq q/2] = \frac{O_s(1)}{(q/2)^s}.\]
We call a sequence of vertices $(w_1, w_2, \dots, w_r)$ bad if there are more than $c_T$ copies of $T$ in $G$ such that $w_i$ corresponds to $u_i$ for all $1 \leq i \leq r$. If we let $B$ be the random variable counting the number of bad sequences, we have, since $s = 2br$ and $q$ is sufficiently large,
\[\mathbb{E}[B] \leq 2 N^r \cdot \frac{O_s(1)}{(q/2)^s} = O_s(q^{br - s}) = o(1).\]
We now remove a vertex from each bad sequence to form a new graph $G'$. Since each vertex has degree at most $N$, the total number of edges removed is at most $B N$. Hence, the expected number of edges in $G'$ is
\[N^{2 - a/b} - \mathbb{E}[B] N = \Omega(N^{2 - a/b}).\]
Therefore, there is a graph with at most $2N$ vertices and $\Omega(N^{2 - a/b})$ edges such that no sequence of $r$ vertices has more than $c_T$ labelled copies of $T$ rooted on these vertices. Finally, we note that this result was only shown to hold when $q$ is a prime power and $N = q^b$. However, an application of Bertrand's postulate shows that the same conclusion holds for all $N$.
\qed

\section{Concluding remarks} \label{sec:conclusion}

We have shown that for any rational number $r$ between $1$ and $2$, there exists a family of graphs $\mathcal{H}_r$ such that $\textrm{ex}(n, \mathcal{H}_r) = \Theta(n^r)$. However, Erd\H{o}s and Simonovits (see, for example,~\cite{E81}) asked whether there exists a single graph $H_r$ such that $\textrm{ex}(n, H_r) = \Theta (n^r)$. Our methods give some hope of a positive solution to this question, but the difficulties now lie with determining accurate upper bounds for the extremal number of certain graphs.

To be more precise, given a rooted tree $(T, R)$, we define $T^p$ to be the graph consisting of the union of $p$ distinct labelled copies of $T$, each of which agree on the set of roots $R$ but are otherwise disjoint. Lemma~\ref{lower} clearly shows that $\textrm{ex}(n, T^p) = \Omega(n^{2 - 1/\rho_T})$ when $T$ is a balanced rooted tree. We believe that a corresponding upper bound should also hold.

\begin{conjecture}
For any balanced rooted tree $(T, R)$, the graph $T^p$ satisfies
\[\textrm{ex}(n, T^p) = O_p(n^{2 - 1/\rho_T}).\]
\end{conjecture}

The condition that $(T, R)$ be balanced is necessary here, as may be seen by considering the graph in Figure~\ref{fig:unbalanced}, namely, a star $K_{1,3}$ with two rooted leaves. Then $T^2$ contains a cycle of length $4$, so the extremal number is $\Omega(n^{3/2})$, whereas the conjecture would suggest that it is $O(n^{4/3})$. 

In order to solve the Erd\H{o}s--Simonovits conjecture, it would be sufficient to solve the conjecture for the collection of rooted trees $T_{a,b}$ with $a < b$ and $(a,b) = 1$. However, even this seems surprisingly difficult and the only known cases are when $a = 1$, in which case $T$ is a star with rooted leaves and $T^p$ is a complete bipartite graph, or $b - a = 1$, when $T$ is a path with rooted endpoints and $T^p$ is a theta graph.

\vspace{3mm}
{\bf Acknowledgements.} We would like to thank Jacques Verstraete for interesting discussions relating to the topic of this paper. We would also like to thank an anonymous referee and Lisa Sauermann for a number of useful comments and corrections.


\begin{thebibliography}{}

\bibitem{ARS99}
{N. Alon, L. R\'onyai and T. Szab\'o,} {Norm-graphs: variations and applications,} {\it J. Combin. Theory Ser. B} {\bf 76} (1999), 280--290.

\bibitem{BBK13}
{P. V. M. Blagojevi\'c, B. Bukh and R. Karasev,} {Tur\'an numbers for $K_{s,t}$-free graphs: topological obstructions and algebraic constructions,} {\it Israel J. Math.} {\bf 197} (2013), 199--214.

\bibitem{B66}
{W. G. Brown,} {On graphs that do not contain a Thomsen graph,} {\it Canad. Math. Bull.} {\bf 9} (1966), 281--285.

\bibitem{B14}
{B. Bukh,} {Random algebraic construction of extremal graphs,} {\it Bull. London Math. Soc.} {\bf 47} (2015), 939--945.

\bibitem{C14}
{D. Conlon,} {Graphs with few paths of prescribed length between any two vertices,} to appear in {\it Bull. London Math. Soc.}

\bibitem{E38}
{P. Erd\H{o}s,} {On sequences of integers no one of which divides the product of two others and on some related problems,} {\it Mitt. Forsch.-Inst. Math. Mech. Univ. Tomsk} {\bf 2} (1938), 74--82.

\bibitem{E81} 
{P. Erd\H{o}s,} {On the combinatorial problems which I would most like to see solved,} {\it Combinatorica} {\bf 1} (1981), 25--42.

\bibitem{ERS66}
{P. Erd\H{o}s, A. R\'enyi and V. T. S\'os,} {On a problem of graph theory,} {\it Studia Sci. Math. Hungar.} {\bf 1} (1966), 215--235.

\bibitem{ES66}
{P. Erd\H{o}s and M. Simonovits,} {A limit theorem in graph theory,} {\it Studia Sci. Math. Hungar.} {\bf 1} (1966), 51--57.

\bibitem{ES46} 
{P. Erd\H{o}s and A. H. Stone,} {On the structure of linear graphs,} {\it Bull. Amer. Math. Soc.} {\bf 52} (1946), 1087--1091.

\bibitem{FS83}
{R. J. Faudree and M. Simonovits,} {On a class of degenerate extremal graph problems,} {\it Combinatorica} {\bf 3} (1983), 83--93.

\bibitem{F86}
{P. Frankl,} {All rationals occur as exponents,} {\it J. Combin. Theory Ser. A} {\bf 42} (1986), 200--206.

\bibitem{FS13}
{Z. F\"uredi and M. Simonovits,} {The history of degenerate (bipartite) extremal graph problems,} in {Erd\H{o}s Centennial,} 169--264, {Bolyai Soc. Math. Stud., 25,} {Springer, Berlin,} 2013.

\bibitem{KRS96}
{J. Koll\'ar, L. R\'onyai and T. Szab\'o,} {Norm-graphs and bipartite Tur\'an numbers,} {\it Combinatorica} {\bf 16} (1996), 399--406.

\bibitem{KST54}
{T. K\H{o}v\'ari, V. T. S\'os and P. Tur\'an,} {On a problem of K. Zarankiewicz,} {\it Colloq. Math.} {\bf 3} (1954), 50--57.

\bibitem{LW54}
{S. Lang and A. Weil,} {Number of points of varieties in finite fields,} {\it Amer. J. Math.} {\bf 76} (1954), 819--827.

\end{thebibliography}
\end{document}